\documentclass[11pt,reqno]{amsart}
\usepackage{graphicx}
\usepackage[hidelinks,pdfusetitle]{hyperref}
\usepackage{amsmath}
\usepackage{amsfonts,amssymb}
\usepackage{amsmath,amssymb,mathtools,mathrsfs}
\usepackage{tikz,longtable}
\usepackage{marvosym}
\newtheorem{theorem}{Theorem}[section]
\newtheorem{lemma}[theorem]{Lemma}
\newtheorem{proposition}[theorem]{Proposition}
\newtheorem{corollary}[theorem]{Corollary}
\theoremstyle{definition}

\theoremstyle{remark}
\newtheorem{remark}[theorem]{Remark}
\numberwithin{equation}{section}

\DeclareRobustCommand{\FieldDegree}{n}
\newcommand{\F}{\mathbb F}

\newcommand{\D}{\mathcal D}
\newcommand{\PSL}{\operatorname{PSL}}

\newcommand{\SU}{\operatorname{SU}}

\newcommand{\SL}{\operatorname{SL}}
\newcommand{\Aut}{\operatorname{Aut}}
\newcommand{\Cos}{\operatorname{Cos}}
\newcommand{\Tr}{\operatorname{Tr}}
\newcommand{\diam}{\operatorname{diam}}

\newcommand{\Spec}{\operatorname{Spec}}
\newcommand{\Ga}{\Gamma(q,\xi)}
\newcommand{\lx}{\ell_\xi}
\newcommand{\px}{\psi_\xi}
\newcommand{\Hx}{H_\xi}

\title[Trace-defined double covers of the  $S_3$-involution graph of $\mathrm{PSL}_2(2^{\FieldDegree})$]{Trace-defined double covers of the  $S_3$-involution graph of $\mathrm{PSL}_2(2^{\FieldDegree})$}
\newcommand{\FirstAuthor}{G. A. Glinskikh}
\newcommand{\SecondAuthor}{L.~Yu.~Tsiovkina}

\author{\FirstAuthor}
\author{\SecondAuthor}
\address[Both authors]{N.N. Krasovskii Institute of Mathematics and Mechanics of UB RAS,\\ 16 S. Kovalevskaya Str., Yekaterinburg,
      620990, Russian Federation}
\address[Both authors]{Ural Federal University,\\51 Lenina Str., Yekaterinburg,
      620000, Russian Federation}
\email[\FirstAuthor]{georgy.glinskikh@imm.uran.ru}
\email[\SecondAuthor]{tsiovkina@imm.uran.ru}
\subjclass[2020]{Primary 05E30; Secondary 05C25}
\keywords{double cover, coset graph, $2$-dimensional special linear group, involution graph,
antipodal distance-regular graph, amply regular graph, Deza graph,  girth,
Ramanujan graph}
\date{}

\begin{document}
\begin{abstract}
For $q=2^n\ge 4$, we construct a family of Deza double covers $\Gamma(q,\xi)$ of the $S_3$-involution graph of $\mathrm{PSL}_2(q)$, parametrized by $\xi\in\mathbb F_q^*$. We show that each resulting cover is connected, arc-transitive and amply regular with parameters $(2(q^2-1),q,\operatorname{Tr}(\xi),1)$, where $\operatorname{Tr}$ denotes the absolute trace. We also show that the trace-one covers are antipodal of diameter $5$ and girth $3$, whereas the trace-zero covers are non-antipodal of diameter $3$ and girth $5$. For $q=4$, these cases yield the well-known icosidodecahedral and dodecadodecahedral graphs, respectively. Then we determine the distance distributions and full automorphism groups of the covers and prove that, for fixed $q$, two covers are isomorphic if and only if their defining parameters lie in the same Frobenius orbit. We also express the Ramanujan property of a cover in terms of the spectral norm of a signed adjacency matrix of the base graph.
\end{abstract}
\maketitle

\section{Introduction}\label{sec:introduction}
\setcounter{equation}{0}
We consider the problem of constructing and classifying coset graphs of finite groups admitting a $BN$-pair of rank $1$, and determining how their geometry depends on the vertex stabiliser and the adjacency relation.
Coset graphs of finite groups admitting a $BN$-pair of rank $1$ provide a natural framework for constructing and studying symmetric covers of complete graphs. The Bruhat decomposition gives a doubly transitive action of such a group $G$ on the right cosets of a Borel subgroup $B$ and, for any subgroup $H\le B$, the natural projection $Hx\mapsto Bx$ gives a $G$-invariant partition.   If one  chooses $H$  as an overgroup of the unipotent radical of $B$, then standard   coset graph construction on $G/H$ gives several infinite families of arc-transitive  antipodal covers of complete graphs \cite{Tsiovkina22}, and describes almost all  examples of distance-regular such covers with  non-solvable automorphism groups. The following  question naturally arises.

\textbf{Problem.}
{\it For which subgroups $H\le B$ one  can choose a $G$-invariant adjacency relation on $G/H$ such that the projection $Hx\mapsto Bx$ defines a cover of the complete graph on $G/B$? What are the combinatorial properties of the resulting graphs?}

In general, covers of complete graphs  are  closely connected with important objects in algebraic combinatorics and discrete geometry, including association schemes and, in the regular abelian case, equiangular line systems (see an overview in \cite{KP2011}). This motivates the study of how the choice of $H$ and the adjacency relation influences the symmetry, distances and regularity of the resulting covers.

The present paper continues the second author's series of studies \cite{Tsi25+,Tsi26+} on coset graphs whose vertex stabiliser $H$ in $G$ has index $2$ in the unipotent radical $U$ of a Borel subgroup $B$, where $G$ is a Suzuki simple group or a unitary group.  With suitable adjacency relations, taking $H$ as an index $2$ subgroup of $U$ produces double  covers of known coset graphs that are antipodal covers of complete graphs. The main questions are how this refinement changes distances and neighbourhood intersections, when antipodal structures arise, and how their full automorphism groups and isomorphism classes can be determined.

Here we treat the linear case $G=\operatorname{PSL}_2(2^n)$. Our construction yields trace-defined double covers of the $S_3$-involution graph of $G$, which belong to the class of Deza graphs. Namely, an absolute trace condition determines two distinct distance geometries: the trace-one graphs form an infinite family of arc-transitive antipodal double covers of diameter $5$, whose smallest member is the icosidodecahedral graph, while the trace-zero graphs are non-antipodal of diameter $3$, whose smallest member is the dodecadodecahedral graph. We determine their local parameters, distance distributions and full automorphism groups, and classify them up to isomorphism. The construction also gives explicit double covers of a Ramanujan quotient, leading to a spectral question about whether the covers themselves retain the Ramanujan property.

Let $q=2^n\ge4$ and $G=\operatorname{PSL}_2(q)$. The $S_3$-involution graph $\Delta_n$ has the involutions of $G$ as its vertices, with two vertices adjacent whenever their product has order $3$. As an abstract graph it was originally discovered by  Mathon (see \cite{BCN}), and its realization as an involution graph is given in \cite[Proposition~4.1 and Remark~4.2]{Tsiovkina2017}; see also \cite{DevillersGiudici}. It is an antipodal distance-regular $(q-1)$-cover of $K_{q+1}$ with intersection array $$\{q,q-2,1;\,1,1,q\},$$ eigenvalues $q,-1$ and $\pm\sqrt{q}$,
and it admits an arc-transitive action of $G$.

Since $q$ is even, we identify $G$ with $\operatorname{SL}_2(q)$ and put
$$u(t)=\begin{pmatrix}1&t\\0&1\end{pmatrix},\qquad S=\{u(t):t\in\mathbb F_q\}, \qquad g=\begin{pmatrix}0&1\\1&0\end{pmatrix}.$$
Here $S$ is the unipotent radical of the Borel subgroup $M=N_G(S)$, and $g$ is a Weyl involution. Denote by $\operatorname{Cos}(G,H,D)$ the graph on the right $H$-cosets of $G$, in which $Hx$ and $Hz$ adjacent whenever $xz^{-1}\in D$, where $D$ is an inverse-closed union of $H$-double cosets disjoint from $H$. In this notation, $$\Delta_n\simeq\operatorname{Cos}(G,S,SgS).$$

Write $\operatorname{Tr}$ for   the absolute trace map $\operatorname{Tr}_{\mathbb F_q/\mathbb F_2}:\mathbb F_q\rightarrow \mathbb F_2$. Since $S$ is elementary abelian, its index-two subgroups are precisely
$$H_\xi=\{u(t):\operatorname{Tr}(\xi t)=0\}, \qquad \xi\in\mathbb F_q^*.$$
For $y\in S- H_\xi$, define
\begin{equation}\label{eq:intro-general-graph}
\Gamma(q,\xi)=\operatorname{Cos}\!\left(G,H_\xi,H_\xi ygH_\xi\cup H_\xi(yg)^{-1}H_\xi\right).
\end{equation}
We show that this graph does not depend on the choice of $y$, and the natural projections $H_\xi x\mapsto Sx\mapsto Mx$ give a tower of covers
\begin{equation}\label{eq:intro-cover-tower}
\Gamma(q,\xi)\longrightarrow\Delta_n\longrightarrow K_{q+1}
\end{equation}
of respective indices $2$ and $q-1$, which implies that $\Gamma(q,\xi)$ is a $q$-regular double cover of $\Delta_n$ on $2(q^2-1)$ vertices, and a $2(q-1)$-cover of $K_{q+1}$. We also show that it admits an arc-regular group of automorphisms $G\times\langle\tau_\xi\rangle$, where
$\tau_\xi:H_\xi x\longmapsto H_\xi yx$
is an involution centralizing $G$ and interchanging the two $G$-orbits on the arcs.

Our first main result determines the local parameters and distance structure of the constructed covers $\Gamma(q,\xi)$.

\begin{theorem}\label{thm:intro-main}
Let $q=2^n\ge4$ and $\xi\in\mathbb F_q^*$. Then $\Gamma(q,\xi)$ is a connected, arc-transitive, amply regular graph with parameters
$$
(v,k,\lambda,\mu)=\bigl(2(q^2-1),q,\operatorname{Tr}(\xi),1\bigr),
$$
and the following statements are true.
\begin{enumerate}
\item[\rm(i)] If $\operatorname{Tr}(\xi)=0$, then $\Gamma(q,\xi)$ has girth $5$, diameter $3$, and distance distribution
$$
\bigl(1,q,q(q-1),q^2-3\bigr).
$$
The two vertices of each covering fibre over $\Delta_n$ are at distance $3$, but $\Gamma(q,\xi)$ is not antipodal.
\item[\rm(ii)] If $\operatorname{Tr}(\xi)=1$, then $\Gamma(q,\xi)$ has girth $3$, every edge lies in a unique triangle, and it is an antipodal double cover of $\Delta_n$ of diameter $5$. Its antipodal classes are the pairs $\{v,\tau_\xi(v)\}$. For $q\ge8$, its distance distribution is
$$
\bigl(1,q,q(q-2),q^2-4,q,1\bigr).
$$
\item[\rm(iii)] For $q=4$, the trace-one graphs are isomorphic to the icosidodecahedral graph,
while the trace-zero graphs are isomorphic to the dodecadodecahedral graph.
\end{enumerate}
\end{theorem}
\begin{corollary}
 For all $q=2^n\ge4$ and $\xi\in\mathbb F_q^*$, $\Gamma(q,\xi)$ is a connected Deza graph with parameters $$(v,k,a,b)=(2(q^2-1),q,0,1).$$
\end{corollary}

Both trace cases occur over every field $\mathbb F_q$ with $q\ge4$. To describe the full automorphism groups, put
$$
d_\xi=[\mathbb F_2(\xi):\mathbb F_2],
\qquad
\mathcal F_\xi=\{\sigma\in\operatorname{Gal}(\mathbb F_q/\mathbb F_2):\sigma(\xi)=\xi\}
\simeq Z_{n/d_\xi}.
$$

\begin{theorem}\label{thm_full_aut}
For each $q=2^n\ge4$ and $\xi\in\mathbb F_q^*$,
\begin{equation}\label{eq:intro-full-aut}
\operatorname{Aut}(\Gamma(q,\xi))\simeq Z_2\times\bigl(\operatorname{PSL}_2(q)\rtimes\mathcal F_\xi\bigr),
\end{equation}
where the central factor $Z_2$ is generated by $\tau_\xi$, and $\mathcal F_\xi$ acts by field automorphisms. In particular,
$$
\begin{aligned}
|\operatorname{Aut}(\Gamma(q,\xi))|&=2q(q^2-1)\frac{n}{d_\xi},\\
|\operatorname{Aut}(\Gamma(q,\xi))_v|&=q\frac{n}{d_\xi},\\
|\operatorname{Aut}(\Gamma(q,\xi))_{(v,w)}|&=\frac{n}{d_\xi}
\end{aligned}
$$
for every vertex $v$ and arc $(v,w)$. In particular, the full automorphism group is arc-regular if and only if $d_\xi=n$.
\end{theorem}

The next theorem shows that the isomorphism classes of resulting covers are determined by the Frobenius action on the defining coefficients.

\begin{theorem}\label{thm_iso}
For each $q=2^n\ge4$ and nonzero $\xi,\eta\in\mathbb F_q$,
\begin{equation}\label{eq:intro-iso-criterion}
\Gamma(q,\xi)\simeq\Gamma(q,\eta)
\quad\Longleftrightarrow\quad
\eta=\xi^{2^j}\text{ for some }j\in\{0,\ldots,n-1\}.
\end{equation}
\end{theorem}

Structural analysis of the constructed covers is concluded with a study of whether they inherit the property of being Ramanujan from their quotients $\Delta_n$.
\begin{theorem}\label{thm_ram} For each $q=2^n\ge4$ and nonzero $\xi\in\mathbb F_q$,
$$
\Gamma(q,\xi)\text{ is Ramanujan} \quad\Longleftrightarrow\quad \|C_\xi\|\le2\sqrt{q-1},
$$
where $C_\xi$ is a signed adjacency matrix  of $\Delta_n$ and  $\|C_\xi\|$ denotes its spectral norm.
\end{theorem}

The paper is organized as follows.
In Section~\ref{sec:preliminaries} we provide  terminology and notation. In  Section~\ref{SecAux}, we  introduce  the construction and prove  its basic properties. Section~\ref{sec:parity} is devoted to computation of  the distance distributions. The proofs combine an analysis of distribution of right $H$-cosets and  paths lifted from $\Delta_n$. Coordinates for the covering graphs lead to trace equations over $\mathbb F_q$, whose solvability is studied using estimates for binary Kloosterman sums. The trace-zero case is treated in Theorem~\ref{thm:even-diameter}, the trace-one case in
Theorem~\ref{thm:odd-diameter}, and the exceptional case $q=4$ in Proposition~\ref{prop:q4-exception}.  In Section~\ref{sec:automorphism} we study automorphisms and isomorphism classes of the constructed  covers proving Theorems~\ref{thm_full_aut} and \ref{thm_iso}. In particular, we give an intrinsic recognition of the covering fibres, including in the non-antipodal case (Lemma~\ref{lem:deck-intrinsic}). In  Section~\ref{sec:deza-girth} we establish ample regularity,
Deza property, the girth formula, and prove Theorem~\ref{thm_ram}.

\section{Notation and terminology}\label{sec:preliminaries}

All graphs in this paper are finite, simple, and undirected. For a graph
$\Gamma$ and vertices $u,v\in V\Gamma$, let $d_\Gamma(u,v)$ denote the
usual graph distance. For $i\ge0$ and $v\in V\Gamma$, put
$$
\Gamma_i(v)=\{x\in V\Gamma:d_\Gamma(v,x)=i\}.
$$
Thus $\Gamma_0(v)=\{v\}$ and $\Gamma_1(v)=\Gamma(v)$ is the neighbourhood
of $v$. We refer to $\Gamma_i(v)$ as the \emph{$i$-th distance
neighbourhood}, or simply the \emph{$i$-th neighbourhood}, around $v$.
The diameter of $\Gamma$ is
$$
\diam(\Gamma)=\max\{d_\Gamma(u,v):u,v\in V\Gamma\}.
$$

A surjective graph morphism
$$
\pi:\widetilde{\Gamma}\longrightarrow\Gamma
$$
is a \emph{covering projection} if, for every
$\widetilde v\in V\widetilde{\Gamma}$, the restriction
$$
\pi:\widetilde{\Gamma}(\widetilde v)
   \longrightarrow \Gamma(\pi(\widetilde v))
$$
is a bijection. The graph $\widetilde{\Gamma}$ is then called a
\emph{cover} of $\Gamma$. For $v\in V\Gamma$, the set
$$
\pi^{-1}(v)
$$
is called the \emph{fibre} over $v$. If every fibre has size $r$, then
$\widetilde{\Gamma}$ is an \emph{$r$-fold cover} of $\Gamma$.

If $v\in V\Gamma$ and $\widetilde v\in\pi^{-1}(v)$, then
$\widetilde v$ is called a \emph{lift} of $v$. More generally, if
$$
W=(v_0,v_1,\ldots,v_m)
$$
is a walk in $\Gamma$ and $\widetilde v_0\in\pi^{-1}(v_0)$, then there
is a unique walk
$$
\widetilde W=(\widetilde v_0,\widetilde v_1,\ldots,\widetilde v_m)
$$
in $\widetilde{\Gamma}$ such that
$$
\pi(\widetilde v_i)=v_i
\qquad (0\le i\le m).
$$
We call $\widetilde W$ the \emph{lift of $W$ starting at
$\widetilde v_0$}.

For a regular cover with covering transformation group $K$, one may
describe the cover by assigning to each oriented edge of the base graph
an element of $K$, called its \emph{voltage}. The voltage of a walk is
the ordered product of the voltages of its oriented edges, with reversed
edges contributing inverse elements. A closed walk has \emph{trivial
voltage} if this product is the identity element of $K$. Equivalently,
a closed walk has trivial voltage if and only if each of its lifts is
closed. For a double cover, where $K\simeq Z_2$, we shall use the terms
\emph{trivial} and \emph{nontrivial voltage} for the two possible
values.

A connected graph $\Gamma$ of diameter $d$ is called
\emph{antipodal} if the relation
$u=v\ \text{ or }\ d_\Gamma(u,v)=d$
is an equivalence relation on $V\Gamma$. Its equivalence classes are
called the \emph{antipodal classes}. If every antipodal class has size
$r$, then $\Gamma$ is said to be \emph{$r$-antipodal}. In particular,
in an antipodal double cover every vertex has a unique vertex at maximum
distance, called its \emph{antipode}.

A connected graph $\Gamma$ of diameter $d$ is called
\emph{distance-regular} if, whenever $u,v\in V\Gamma$ satisfy
$d_\Gamma(u,v)=i$, the numbers
$$c_i=|\Gamma(v)\cap\Gamma_{i-1}(u)|,\qquad a_i=|\Gamma(v)\cap\Gamma_i(u)|,\qquad b_i= |\Gamma(v)\cap\Gamma_{i+1}(u)|$$
depend only on $i$, and not on the particular choice of $u$ and $v$.
Here
$$c_0=0,\qquad b_d=0.$$
The numbers $b_i,c_i$ are the \emph{intersection numbers}, and the
sequence
$$\{b_0,b_1,\ldots,b_{d-1};c_1,c_2,\ldots,c_d\}$$
is called the \emph{intersection array} of $\Gamma$. If $\Gamma$ has
valency $k=b_0$, then
$$a_i=k-b_i-c_i.$$
A regular graph on $v$ vertices of valency $k$ is a Deza graph with
 parameters $(v,k,b,a)$ if every two distinct vertices have either $b$ or
 $a$ common neighbours, where $b\ge a$.; see \cite{Deza99}.
 We write $\operatorname{girth}(\Gamma)$ to denote the girth of $\Gamma$, i.e., the length of its shortest
 cycle.

For a group $G$, its core-free subgroup $H$ and an inverse-closed
union $\D$ of $H$-double cosets with $1\notin \D$, we define the coset graph $\Cos(G,H,\D)$ on
 the right-coset space $G/H$ (where the right $H$-cosets are written as $Hx$)
 by the following rule $$ Hx\sim Hz\quad\Longleftrightarrow\quad xz^{-1}\in \D.$$
The group $G$ acts (faithfully)   by right multiplication, since $(xr)(zr)^{-1}=xz^{-1}$. As $\D=\D^{-1}$,
it follows that the neighbours of $Hx$ can   be written as $Hdx$ with $d\in D$. For more background on characterization of vertex-transitive graphs as coset graphs, see
\cite{BeinekeWilson}.
\section{The construction and auxiliary results}\label{SecAux}
Throughout the rest we let $q=2^n\ge4$ and put
$$G=\SL_2(q),\qquad
u(t)=\begin{pmatrix}1&t\\0&1\end{pmatrix},\quad
S=\{u(t):t\in\F_q\},\quad
g=\begin{pmatrix}0&1\\1&0\end{pmatrix}.$$

Since $q$ is even,  $G$ has trivial center, so it can be also identified with $\PSL_2(q)$.
Let $K$ be the diagonal torus and $M=S\rtimes K$ the Borel subgroup  of $G$.
The monomial involution $g\in G-S$ represents the nontrivial Weyl element.
We use $x$ for a general element of $G$ and write $u(t)$ for elements of $S$.
Let
$$\Tr=\Tr_{\F_q/\F_2}:\F_q\longrightarrow \F_2$$
denote the absolute trace map,
$$\Tr(x)=x+x^2+x^{2^2}+\cdots+x^{2^{n-1}}.$$
Put $\psi(t)=(-1)^{\Tr(t)}$.
For $\xi\in \mathbb{F}_q^*$, we define
$$ \lx(t)=\Tr(\xi t),\qquad E_\xi=\ker\lx,\qquad \Hx=\{u(t):t\in E_\xi\},$$
$$\px(t)=(-1)^{\lx(t)},\qquad \delta_\xi=\lx(1)=\Tr(\xi).$$

The maps $\ell_\xi$, where $\xi\in\F_q$, are exactly the $\F_2$-linear maps
$\F_q\to\F_2$. Since distinct nonzero $l_{\xi}$ have distinct kernels, $\{H_\xi:\xi\ne0\}$ forms the set of all index-two subgroups of $S$.
There are $q/2$ coefficients with $\delta_\xi=1$ and $q/2-1$ nonzero coefficients with $\delta_\xi=0$.
\medskip

Choose $c\in\F_q$ with $\lx(c)=1$ and set $$y=u(c)\in S-H_\xi,\qquad a_c=yg.$$
Define
\begin{equation}\label{eq:general-graph}
\Ga=\Cos\bigl(G,H_\xi, H_\xi \{yg, (yg)^{-1}\}H_\xi\bigr).
\end{equation}

\begin{lemma}\label{lem:c-independent}
For a fixed $\xi\in\F_q^*$, the graph in~\eqref{eq:general-graph} is independent on the choice of
$c\in\F_q^*$ such that $\ell_\xi(c)=1$. The map
$$\tau_\xi:H_\xi x\longmapsto H_\xi yx$$
is also independent of this choice.
\end{lemma}
\begin{proof}
If $\ell_\xi(c')=1$, then $c'=c+e$ with $u(e)\in H_\xi$.
Thus $a_{c'}=u(e)a_c$ and $a_{c'}^{-1}=a_c^{-1}u(e)$, which do not change
the two double cosets. The same equality gives
$H_\xi u(c')x=H_\xi u(c)x$.
\end{proof}

\begin{proposition}\label{prop:basic}
For every $\xi\in\F_q^*$, the graph $\Ga$ has $2(q^2-1)$ vertices and valency
$q$. The natural projection maps
$$G/H_\xi\longrightarrow G/S \longrightarrow G/M,$$
give a sequence of covers
$$\Ga\longrightarrow\Delta_n\longrightarrow K_{q+1},$$
where the first graph quotient is $\Delta_n=\Cos(G,S,SgS)$ and the second
is $K_{q+1}=\Cos(G,M,MgM)$. The respective covering indices are   $2$ and $q-1$, hence
$\Ga$ is an $2(q-1)$-cover of $K_{q+1}$.
\end{proposition}
\begin{proof}
By assumption we have $|G:H_\xi|=2(q^2-1)$. Define
$$b(t)=c(1+\lx(t)),\qquad d(t)=u(b(t))gu(t)  =\begin{pmatrix}b(t)&1+b(t)t\\1&t\end{pmatrix}.$$
For $t\in E_\xi$, we have  $d(t)=a_cu(t)$, while  for $t\notin E_\xi$,
we have $d(t)=a_c^{-1}u(t+c)$. Thus the $q$ neighbours of $H_\xi$ are exactly
$H_\xi d(t)$, where $t\in\F_q$. They are distinct because their lower rows
$(1,t)$ are distinct. Their images are the $q$ distinct neighbours
$Sg u(t)$ of $S$ in $\Delta_n$. The projective lower rows $[1:t]$ are
exactly the projective points other than $[0:1]$. This shows local
bijectivity for both projection maps with respect to the corresponding vertices $H_\xi$ and $S$,
which by transitivity of  the right action of $G$ on $G/H_\xi$ implies that both maps are locally
bijective with respect to every vertex.
 \end{proof}

\begin{proposition}\label{prop:arc-transitive}
$G$ acts vertex- and edge-transitively, with two orbits on arcs, and the fixed-point-free involution $\tau_\xi$ is an automorphism of $\Ga$, which centralizes $G$  and interchanges the two $G$-orbits on arcs.
\end{proposition}
\begin{proof}
If $u(r)u(c)gu(s)=gu(c)$ for $r,s\in E_\xi$, comparison of the
$(1,1)$-entries gives $r=c$, a contradiction. Also
$H_\xi\cap H_\xi^{a_c}=1$, so each double coset determines a directed
suborbit of size $q/2$. By definition of $\Ga$, the directed
orbital indexed by $H_\xi a_cH_\xi$ contains $(H_\xi,H_\xi a_c^{-1})$,
whereas $(H_\xi,H_\xi a_c)$ belongs to the inverse orbital.
Since $S$ is elementary abelian, $y=u(c)$ normalizes $H_\xi$, and $y$ is an involution
not contained in $H_\xi$. Note also that conjugation by $y$ interchanges the two  double cosets in $D$.
Hence  for the images of $H_\xi x,H_\xi z$ under $\tau_\xi$, we have
$$(yx)(yz)^{-1}=y(xz^{-1})y^{-1}\in \D^y=\D$$ whenever $xz^{-1}\in \D$.
So the involution $\tau_\xi$ is fixed-point-free. Moreover, it centralizes $G$.
Since  $G$ is simple, it acts faithfully on $G/M$, and we have  $\tau_\xi\not \in G$. The rest follows.
\end{proof}

We will say that the vertices $H_\xi x,H_\xi z$ are {\it $\tau_\xi$-opposite}, if they belong to a common $\langle\tau_\xi\rangle$-orbit.
\begin{remark}\label{rem:fixed-y}
In what follows we  take $c=1$ whenever $\Tr(\xi)=1$. However, if $\xi\ne0$ but
$\Tr(\xi)=0$, this is not an allowed choice as defining double cosetes coincide in this case.

 \begin{proposition}
If $g_1=h(r)g$ with $r\in\F_q^*$, then
 $$\Cos(G,H,H\{g_1y,(g_1y)^{-1}\}H)\simeq \Gamma(q,\xi r).$$
 Hence  the family of graphs $\Ga$ is independent (up to isomorphism) on the choice of involution  $g\in G-S$.
 \end{proposition}
\begin{proof}
Since $q$ is even, the squaring map is an automorphism of
$\F_q^*$. Choose $s\in\F_q^*$ such that $s^2=r$.
Then
$$g_1=h(s)gh(s)^{-1}.$$
Consider the inner automorphism
$$\varphi_s:x\longmapsto h(s)^{-1}xh(s).$$
It satisfies
$$\varphi_s(g_1)=g$$
and
$$\varphi_s(u(t))=u(s^{-2}t).$$
Hence
$$\varphi_s(H_\xi)=H_{\xi s^2}=H_{\xi r}.$$

Put $c'=s^{-2}c$. Since $\ell_\xi(c)=1$, we have
$$\ell_{\xi r}(c')=\Tr(\xi r s^{-2}c)=\Tr(\xi c) =1.$$
Thus $\varphi_s(y)=u(c')$ is an admissible choice of the element
defining $\Gamma(q,\xi r)$.

Moreover,
$$\varphi_s\bigl(\{g_1y,(g_1y)^{-1}\}\bigr)=\{g u(c'),u(c')g\}=\{u(c')g,(u(c')g)^{-1}\}.$$
Therefore $\varphi_s$ maps the defining double cosets onto those of
$\Gamma(q,\xi r)$ and induces a graph isomorphism
$$\Cos\bigl(G,H_\xi, H_\xi\{g_1y,(g_1y)^{-1}\}H_\xi\bigr) \simeq\Gamma(q,\xi r).$$

Finally, multiplication by $r$ permutes $\F_q^*$. Hence replacing
$g$ by $g_1$ essentially  reparametrizes the family
$\{\Gamma(q,\xi):\xi\in\F_q^*\}$.
\end{proof}

\end{remark}

\section{Distance distribution of $\Ga$}\label{sec:parity}
In three subsections below we assume $q\ge8$.
The case $q=4$ is treated at the end.

\subsection{Coordinates and the common base graph}
Observe that a right coset of $S$ is determined by the nonzero lower row  of a
matrix in $G$, and the two $H_\xi$-cosets with fixed lower row $(C,D)$ are distinguished by the next binary invariant:
\begin{equation}\label{eq:coset-invariant}
\epsilon(x)=
\begin{cases}
\ell_\xi(A/C),&C\ne0,\\
\ell_\xi(B/D),&C=0,
\end{cases}
\qquad
 x=\begin{pmatrix}A&B\\C&D\end{pmatrix}.
\end{equation}
Indeed, the upper rows may differ by a multiple of the lower row, and
left multiplication by $H_\xi$ changes that multiple by an element of
$E_\xi$. The involution $\tau_\xi$ maps $\epsilon$ to $\epsilon+1$.

\begin{proposition}\label{thm:involution}
For $q\ge4$, the quotient graph $\Delta_n$ of $\Ga$ has vertex set $G/S$, and
its adjacency relation is given by $$Sx\sim Sz \iff xz^{-1}\in SgS.$$ It is
distance-regular with intersection array $$\{q,q-2,1;1,1,q\},$$
with distance distribution  $1,q,q(q-2),q-2$.
\end{proposition}
\begin{proof}
See \cite[Proposition~4.1]{Tsiovkina2017}.
 \end{proof}

Note also that under $Sx\mapsto u(1)^x$ this is the involution graph of \cite[Proposition~4.1]{Tsiovkina2017}, as the double coset $SgS$ corresponds
to pairs whose product has order $3$. The parameters of $\Delta_n$ imply  that every edge lies in a unique triangle and that there are no
$4$-cycles, hence every pair of distinct vertices of $\Delta_n$  has at most one common neighbour.

\begin{lemma}\label{lem:two-step-param}
Every lower row $(C,D)$ with $C\ne0$ occurs for exactly one word
$d(s)d(t)$. Its parameters are
$$t=\frac{D+1}{C},\qquad s=C+b(t).$$

For $C=1$, there are exactly two $H_\xi$-cosets with lower row
$(1,D)$. Exactly one of them is adjacent to $H_\xi$, namely
$H_\xi d(D)$, and its binary invariant is $$1+\ell_\xi(D).$$
Its $\tau_\xi$-opposite has invariant $$\ell_\xi(D).$$
On the other hand, the unique two-step word $d(s)d(t)$ with lower row
$(1,D)$ has binary invariant
$$\ell_\xi(D)+\delta_\xi.$$

Consequently, if $\delta_\xi=0$, the two-step word reaches the
$\tau_\xi$-opposite of the adjacent lift, whereas if $\delta_\xi=1$,
it reaches the adjacent lift itself.
\end{lemma}
\begin{proof}
A direct multiplication gives
$$d(s)d(t)=
\begin{pmatrix}
1+b(s)C & *\\
C & 1+tC
\end{pmatrix},
\qquad C=s+b(t).$$
Hence, for a prescribed lower row $(C,D)$ with $C\ne0$,
$$D=1+tC,$$
and therefore
$$t=\frac{D+1}{C},\qquad s=C+b(t).$$
Thus the two-step word is uniquely determined.

Now let $C=1$. Then the binary invariant of the two-step representative is
$$\ell_\xi\bigl(1+b(s)\bigr)=\delta_\xi+\ell_\xi(b(s)).$$
Since
$$\ell_\xi(b(v))=1+\ell_\xi(v),$$
we get
$$\ell_\xi\bigl(1+b(s)\bigr)=\delta_\xi+1+\ell_\xi(s).$$
But $s=1+b(t)$ and $t=D+1$, so
$$\ell_\xi(s)=\delta_\xi+\ell_\xi(b(t))=\delta_\xi+1+\ell_\xi(t)=\delta_\xi+1+\ell_\xi(D)+\delta_\xi=1+\ell_\xi(D).$$
Hence the two-step invariant equals
$$\ell_\xi(D)+\delta_\xi.$$

The adjacent lift with lower row $(1,D)$ is $H_\xi d(D)$, whose
upper-left entry is $b(D)$. Therefore its binary invariant is
$$\ell_\xi(b(D))=1+\ell_\xi(D),$$
while its $\tau_\xi$-opposite has the complemmatary invariant
$$\ell_\xi(D).$$
This proves the lemma.
 \end{proof}

\begin{lemma}\label{lem:three-step-kloosterman}
Let $C\ne0,1$ and $D\in\F_q$. For $L\ne0$ set
$$t=(L+D)/C,\quad s=(C+1)/L+b(t),\quad r=L+b(s).$$
Then $d(r)d(s)d(t)$ has lower row $(C,D)$ and binary invariant
\begin{equation}\label{eq:general-three}
\epsilon_L=1+\ell_\xi(D/C)+\ell_\xi\!\left(\frac{C+1}{C}L+\frac{(C+1)^2}{CL}\right).
\end{equation}
Both possible binary invariants occur. For $C=0$ and $D\ne0$, the
three-step words with that lower row have invariants
\begin{equation}\label{eq:general-base}
\epsilon_0=1+\ell_\xi(D)+\ell_\xi\bigl(t(1+D^{-2})\bigr),\qquad t\in\F_q.
\end{equation}
In particular, both lifts over $(0,D)$ with $D\ne1$ occur in three steps.
\end{lemma}

\begin{proof}
Put
$$L=r+b(s),\qquad C'=s+b(t).$$
A direct multiplication gives
$$d(r)d(s)d(t)=
\begin{pmatrix}
A&B\\ C&D
\end{pmatrix},$$
where $$C=1+C'L,\qquad D=L+tC,$$
and $$A=C'+b(r)C,\qquad B=(1+b(r)L)(1+tC')+b(r)t.$$

Suppose first that $C\ne0,1$. With
$$C'=\frac{C+1}{L}, \qquad t=\frac{L+D}{C}, $$
the binary invariant is
$$\epsilon_L=\ell_\xi(A/C)=\ell_\xi\!\left(\frac{C+1}{CL}\right)+\ell_\xi(b(r)).$$
Using $$\ell_\xi(b(v))=1+\ell_\xi(v), $$
together with $$ r=L+b(s), \qquad s=\frac{C+1}{L}+b(t), $$
we obtain
$$\begin{aligned}
\epsilon_L&= 1+\ell_\xi(D/C) +\ell_\xi\!\left(L+\frac{L}{C}+\frac{C+1}{L}+\frac{C+1}{CL}\right)\\
&=1+\ell_\xi(D/C)+\ell_\xi\!\left(\frac{C+1}{C}L+\frac{(C+1)^2}{CL}\right),
\end{aligned}$$
which is \eqref{eq:general-three}.

Writing $\psi(x)=(-1)^{\Tr(x)}$, the difference between the numbers
of $L\in\F_q^*$ for which $\epsilon_L=0$ and $\epsilon_L=1$
is, up to sign, $$\sum_{L\in\F_q^*}\psi\!\left(\xi\frac{C+1}{C}L+\xi\frac{(C+1)^2}{C}L^{-1}\right).$$
Since $\xi\ne0$ and $C\ne0,1$, both coefficients
$$\xi\frac{C+1}{C},\qquad\xi\frac{(C+1)^2}{C}$$
are nonzero. Hence this is a non-degenerate binary Kloosterman sum,
whose absolute value is at most $2\sqrt q$; see
\cite{ConradKloosterman}. Since $$q-1>2\sqrt q$$
for $q\ge8$, neither of the two solution counts can vanish. Thus
both binary invariants occur.

Now let $C=0$ and $D\ne0$. From
$$C=1+C'L,\qquad D=L+tC$$
we obtain
$$L=D,\qquad C'=D^{-1}.$$
Using the formula for $B$ gives
$$\frac{B}{D}=D^{-1}+tD^{-2}+b(r).$$
Moreover,
$$r=D+b(s),\qquad s=D^{-1}+b(t).$$
Therefore
$$
\begin{aligned}
\ell_\xi(b(r))&=1+\ell_\xi(r)\\
&=1+\ell_\xi(D)+\ell_\xi(b(s))\\
&=\ell_\xi(D)+\ell_\xi(s)\\
&=\ell_\xi(D)+\ell_\xi(D^{-1})+1+\ell_\xi(t).
\end{aligned}
$$
Hence
$$
\begin{aligned}
\epsilon_0=\ell_\xi(B/D)
&=\ell_\xi(D^{-1})+\ell_\xi(tD^{-2})+\ell_\xi(b(r))\\
&=1+\ell_\xi(D)+\ell_\xi\bigl(t(1+D^{-2})\bigr),
\end{aligned}
$$
which proves \eqref{eq:general-base}.

If $D\ne1$, then
$$1+D^{-2}\ne0.$$
Since $\xi\ne0$, the map
$$t\longmapsto\ell_\xi\bigl(t(1+D^{-2})\bigr)=\Tr\!\bigl(\xi(1+D^{-2})t\bigr)$$
is a nonzero $\F_2$-linear functional on $\F_q$, and therefore takes
both values $0$ and $1$. Hence both lifts over $(0,D)$ occur among
three-step words.
 \end{proof}

\subsection{The trace-zero case}
\begin{theorem}\label{thm:even-diameter}
Let $q\ge8$ and let $\xi\ne0$ satisfy $\Tr(\xi)=0$. Then $\Ga$ is
connected, has diameter $3$, and has distance distribution
$$1,\quad q,\quad q(q-1),\quad q^2-3.$$
For every vertex $v$, its $\tau_\xi$-opposite $\tau_\xi(v)$ is at
distance $3$ from $v$. However, the pairs
$\{v,\tau_\xi(v)\}$ are not antipodal classes, and neither are the
fibres of the covering onto $K_{q+1}$.
\end{theorem}

\begin{proof}
Here $\delta_\xi=\ell_\xi(1)=0$.

For $C=1$, Lemma~\ref{lem:two-step-param} shows that the unique
two-step lift is the $\tau_\xi$-opposite of the lift adjacent to
$H_\xi$. Thus the $q$ rows $(1,D)$ contribute $q$ vertices at
distance exactly $2$.

For $C\ne0,1$, Lemma~\ref{lem:two-step-param} gives exactly one
two-step lift over each lower row $(C,D)$, while
Lemma~\ref{lem:three-step-kloosterman} shows that both lifts are
reached within three steps. Hence the second lift lies at distance
at most $3$.

For $C=0$ and $D\ne1$,
Lemma~\ref{lem:three-step-kloosterman} shows that both lifts over
$(0,D)$ occur in three steps. Finally, for $(C,D)=(0,1)$,
\eqref{eq:general-base} gives
$$\epsilon_0=1+\delta_\xi=1,$$
so the nontrivial lift over $(0,1)$, namely $\tau_\xi(H_\xi)$, is
also reached in three steps.

Any vertex whose lower row is projectively dependent with that of
$H_\xi$ lies in the same fibre of the covering onto $K_{q+1}$.
Distinct vertices in one fibre are neither adjacent nor have a
common neighbour, since between any two fibres the covering induces
a perfect matching. Hence every nontrivial such vertex is at
distance at least $3$ from $H_\xi$. In particular,
$$d(H_\xi,\tau_\xi(H_\xi))=3,$$
and, since $\tau_\xi$ commutes with the transitive action of $G$,
$$d(v,\tau_\xi(v))=3$$
for every vertex $v$.

It follows from the preceding paragraphs that every vertex is at
distance at most $3$ from $H_\xi$, while $\tau_\xi(H_\xi)$ is at
distance exactly $3$. Thus
$$\diam(\Ga)=3,$$
and the graph is connected.

There are $q$ vertices at distance $1$. For every lower row
$(C,D)$ with $C\ne0$, exactly one lift is obtained by a two-step
word. There are
$$q(q-1)$$
such lower rows, and in the case $C=1$ these two-step lifts are the
$\tau_\xi$-opposites of the adjacent lifts, hence are disjoint from
the neigbours of $H_\xi$. Therefore
$$|\Ga_2(H_\xi)|=q(q-1).$$
Since
$$|V\Ga|=2(q^2-1),$$
the remaining vertices form the third neighbourhood of $H_\xi$, of size
$$2(q^2-1)-1-q-q(q-1)=q^2-3.$$
Thus the distance distribution is
$$1,\quad q,\quad q(q-1),\quad q^2-3.$$

The pairs $\{v,\tau_\xi(v)\}$ cannot be antipodal classes, since
each vertex has $q^2-3>1$ vertices at maximum distance. In fact,
if the distance-$3$ relation were an antipodal equivalence
relation, every antipodal class would have size
$$1+|\Ga_3(v)|=q^2-2.$$
But
$$q^2-2\nmid 2(q^2-1),$$
so this is impossible. Hence $\Ga$ is not antipodal. In particular,
the fibres of the covering onto $K_{q+1}$ are not antipodal classes.
 \end{proof}

\subsection{The trace-one case}
\begin{theorem}\label{thm:odd-diameter}
Let $q\ge8$ and $\Tr(\xi)=1$. Then $\Ga$ is a connected antipodal graph of diameter $5$,
with antipodal classes precisely the $\tau_{\xi}$-opposite pairs
$$\{H_\xi x,H_\xi u(1)x\}.$$
Its distance distribution is
$$1,\quad q,\quad q(q-2),\quad q^2-4,\quad q,\quad1.$$
\end{theorem}
\begin{proof}
Here
$$\delta_\xi=\ell_\xi(1)=1.$$
By Lemma~\ref{lem:c-independent}, we may take $c=1$. Thus
$$y=u(1),\qquad a=u(1)g,$$
and a direct calculation gives $a^3=1$.

Consider the double covering
$$\rho:\Ga\longrightarrow\Delta_n.$$
Every triangle of $\Delta_n$ has trivial voltage. Indeed, the triangle
$$S,\ Sa,\ Sa^2$$
lifts through $H_\xi$ to
$$H_\xi,\ H_\xi a,\ H_\xi a^2,\ H_\xi a^3=H_\xi,$$
and, since $G$ is transitive on the arcs of $\Delta_n$ and every edge
lies in a unique triangle, all triangles are $G$-conjugate.

Moreover, $\Delta_n$ has no $4$-cycles. Hence every closed walk of
length at most $4$ has trivial voltage: closed $3$-walks are triangles,
while a closed $4$-walk contains an immediate reversal and reduces,
after deleting backtracking pairs, to the trivial walk. Consequently,
no path of length at most $4$ can join the two lifts in one
$\rho$-fibre. In particular,
$$d(H_\xi,\tau_\xi(H_\xi))\ge5.$$
Similarly, if $v$ is adjacent to $H_\xi$, then its $\tau_\xi$-opposite
cannot be reached from $H_\xi$ in at most three steps, since adjoining
the edge $H_\xi v$ would produce a closed base walk of length at most
$4$ with nontrivial voltage. Hence $$d(H_\xi,\tau_\xi(v))\ge4.$$

We now determine the first three nontrivial  $i$-neighbourhoods of $H_{\xi}$.
For every lower row $(C,D)$ with $C\ne0,1$,
Lemma~\ref{lem:two-step-param} gives one lift at distance $2$, while
Lemma~\ref{lem:three-step-kloosterman} shows that the other lift is
reached in three steps. Thus these rows contribute $q(q-2)$
vertices to each of the second and third neighbourhoods.

For $C=0$ and $D\ne0,1$, both lifts over $(0,D)$ occur in three steps,
contributing a further $2(q-2)$
vertices to the third neighbourhood.

Finally, when $C=1$, Lemma~\ref{lem:two-step-param} and $\delta_\xi=1$ show that the unique two-step word reaches the same lift
which is already adjacent to $H_\xi$. Hence no new vertex from the second neighbourhood arises from these $q$ lower rows,
while their $\tau_\xi$-opposite lifts have not yet been reached.

Therefore $$|\Ga_1(H_\xi)|=q,\qquad|\Ga_2(H_\xi)|=q(q-2),$$
and $$|\Ga_3(H_\xi)|=q(q-2)+2(q-2)=q^2-4.$$

The vertices not yet reached are precisely
$$\tau_\xi(H_\xi)\quad\text{and}\quad W:=\tau_\xi(\Ga(H_\xi)).$$
Indeed, $W$ consists of the $\tau_\xi$-opposite lifts over the $q$
lower rows $(1,D)$, and therefore $|W|=q$.

Let $v\in W$. Since $\tau_\xi$ is an automorphism,
$v$ is adjacent to $\tau_\xi(H_\xi)$. Moreover, the local graph at
$\tau_\xi(H_\xi)$ is a disjoint union of $q/2$ edges, so $v$ has
exactly one neighbour in $W$, namely the third vertex of the unique
triangle containing the edge
$$v\,\tau_\xi(H_\xi).$$
The remaining $q-2$ neighbours of $v$ lie outside
$$W\cup\{\tau_\xi(H_\xi)\}.$$
All such vertices have already been reached in $i$-neighbourhoods with $i=0,1,2,3$.
Since $d(H_\xi,v)\ge4$, any neighbour of $v$ has distance at least $3$
from $H_\xi$. Hence these $q-2$ neighbours lie in $\Ga_3(H_\xi)$.
It follows that
$$d(H_\xi,v)=4$$
for every $v\in W$.

Since $\tau_\xi(H_\xi)$ is adjacent to every vertex of $W$, we obtain
$$d(H_\xi,\tau_\xi(H_\xi))\le5.$$
Together with the earlier lower bound this gives
$$d(H_\xi,\tau_\xi(H_\xi))=5.$$

Thus all vertices have been reached, so $\Ga$ is connected and $$\diam(\Ga)=5.$$
The distance distribution relative to $H_\xi$ is
$$1,\quad q,\quad q(q-2),\quad q^2-4,\quad q,\quad1.$$
By vertex-transitivity the same holds at every vertex. Since the unique
vertex at distance $5$ from $H_\xi x$ is $$\tau_\xi(H_\xi x)=H_\xi u(1)x,$$
the antipodal classes are precisely $$\{H_\xi x,H_\xi u(1)x\}.$$
 \end{proof}

\begin{corollary}\label{lem:trace}
For $\Gamma_n:=\Gamma(2^n,1)$ and $q\ge8$, the two preceding theorems
specialize to diameter $3$ for even $n$ and to the antipodal diameter-$5$
case for odd $n$.
\end{corollary}
\begin{proof}
Here $\delta_1=\Tr(1)=n\pmod2$.
 \end{proof}

\begin{proposition}\label{prop:q4-exception}
For $q=4$, all three nonzero coefficients give connected $4$-regular
double covers of $\Delta_2$. The trace-zero coefficient $\xi=1$ gives
diameter $3$ and distribution $(1,4,12,13)$. The two trace-one coefficients
give antipodal diameter-$5$ graphs with distribution
$(1,4,8,8,8,1).$
\end{proposition}
   \begin{proof} Computation in GAP \cite{GAP}.

  \end{proof}

 \section{Automorphism groups and isomorphism classes}
\label{sec:automorphism}

In this section, we determine the full automorphism groups of the covers $\Ga$  and classify the covers for fixed $q$ up to abstract
isomorphism.

Put
$$ d_\xi=[\F_2(\xi):\F_2],\qquad \mathcal F_\xi=  \{\sigma\in\operatorname{Gal}(\F_q/\F_2):\sigma(\xi)=\xi\}\simeq Z_{n/d_\xi}.$$
We use the known description
\begin{equation}\label{eq:aut-quotient}
 \Aut(\Delta_n)\simeq\operatorname{P\Gamma L}_2(q) =G\rtimes\operatorname{Gal}(\F_q/\F_2),
\end{equation}
of the automorphism group of the $S_3$-involution graph; e.g. see~\cite[Theorem~3.11]{DevillersGiudici}.

\begin{lemma}\label{lem_autbasic}
For the right action of $G$ on $G/H_\xi$,
$$ N_G(H_\xi)=S,\qquad C_{\operatorname{Sym}(V\Gamma(q,\xi))}(G) \simeq N_G(H_\xi)/H_\xi =\langle\tau_\xi\rangle\simeq Z_2. $$
\end{lemma}
\begin{proof}
Every non-identity element of $H_\xi$ fixes exactly one point of $\operatorname{PG}(1,q)$, namely the point stabilized by the upper
triangular Borel group $M=S\rtimes K$. Thus a normalizer of $H_\xi$ is contained in $M$.  If $h(s)=\operatorname{diag}(s,s^{-1})\in K$, then $h(s)^{-1}u(t)h(s)=u(s^{-2}t)$. Hence $h(s)$ normalizes $H_\xi$ if and only if $$\ker\ell_{\xi s^{-2}}=\ker\ell_\xi.$$
Nonzero $\F_2$-linear functionals with equal kernels are equal, so $\xi s^{-2}=\xi$ and therefore $s=1$. Conversely, $S$ normalizes
$H_\xi$, since it is abelian.  Hence $N_G(H_\xi)=S$. Finally, the centralizer formula for a transitive coset action gives the
second assertion, with the nontrivial coset in $S/H_\xi$ represented by $\tau_\xi$.
 \end{proof}

\begin{lemma}
\label{lem:deck-intrinsic}
Every automorphism of $\Gamma(q,\xi)$ preserves the partition into pairs $\{v,\tau_\xi(v)\}$, and $\tau_\xi\in Z(\Aut(\Gamma(q,\xi)))$.
If $\Tr(\xi)=0$, then $\tau_{\xi}$-opposite  pairs admit the following  description: for any distinct nonadjacent vertices $v,w$,
\begin{equation}\label{eq_walks}
 (A_{\Gamma(q,\xi)}^3)_{v,w}=q \Longleftrightarrow  w=\tau_\xi(v).
\end{equation}
\end{lemma}
\begin{proof}
For $\Tr(\xi)=1$, the two vertices of every $\tau_{\xi}$-opposite   pair are precisely antipodes (also for the exceptional case $q=4$), and therefore the
partition is invariant under all automorphisms of $\Gamma(q,\xi)$.

Next we assume $\Tr(\xi)=0$.  By vertex-transitivity it is enough to prove~\eqref{eq_walks} for $v=H_\xi$.
As in the proof of Proposition~\ref{prop:basic}, put
$$ b(t)=c(1+\ell_\xi(t)),\qquad  d(t)=u(b(t))gu(t)=
 \begin{pmatrix}
 b(t)&1+b(t)t\\1&t
 \end{pmatrix}.$$
Note that walks of length three from $H_\xi$ correspond bijectively to triples $(r,s,t)\in\F_q^3$ and has an end vertex of the form $H_\xi d(r)d(s)d(t)$.
Indeed,   for every $x\in G$, the set $$\{H_\xi d(t)x\mid t\in \F_q\}$$ is $1$-neighbourhood of $Hx$, and every edge $$H_\xi\sim H_\xi d(s)$$ is mapped to the edge $$H_\xi d(t)\sim H_\xi d(s)d(t)$$ by $d(t)$,
and every edge  $$H_\xi\sim H_\xi d(r)$$ is mapped to the edge $$Hd(s)d(t)\sim H_\xi d(r)d(s)d(t)$$ by $d(s)d(t)$.

Write $(C,D)$ for the lower row of the product $d(r)d(s)d(t)$. Then proof of Lemma~\ref{lem:three-step-kloosterman} gives
$$ C=1+C'L,\quad D=L+tC,\qquad L=r+b(s),\quad C'=s+b(t).$$
We now separate four types of lower rows.

If $C\notin\{0,1\}$, the parameter $L\ne0$ determines a unique triple $(r,s,t)$ for any prescribed $(C,D)$. Thus either of the two
lifts with this lower row receives at most $q-1$ walks of length three. If $C=0$ and $D\ne1$, then $D\ne0$ and formula \eqref{eq:general-base} shows that,
as $t$ varies, the lift coordinate is
$$ \epsilon_0=1+\ell_\xi(D)+\ell_\xi\bigl(t(1+D^{-2})\bigr).$$
As $1+D^{-2}\ne0$, it follows that exactly $q/2$ values of $t$ reach each lift. For the lower row $(C,D)=(0,1)$, this expression is constant,
$1+\ell_\xi(1)=1$. The $q$ corresponding triples all reach $\tau_\xi(H_\xi)$, and none reaches $H_\xi$.

Finally, if $C=1$, then only one lift is adjacent to $H_\xi$, namely $H_\xi d(D)$. The other
lift cannot be an end vertex of a walk of length three, since otherwise, joining the corresponding base vertices in $\Delta_n$ by an edge would give a closed walk of length four with nontrivial voltage.
Every closed walk of length four in $\Delta_n$ backtracks, since $\Delta_n$ has no four-cycles, and consequently has trivial voltage.
It follows that among nonadjacent vertices only its $\tau_{\xi}$-opposite  vertex is the end vertex of precisely $q$  walks of length three. This proves
\eqref{eq_walks}.

In either case, the   partition into $\tau_{\xi}$-opposite pairs is therefore is invariant with respect to every automorphism of  $\Gamma(q,\xi)$. Moreover, every automorphism preserving all $\tau_{\xi}$-opposite pairs commutes with  $\tau_\xi$.
 \end{proof}

\begin{lemma}
\label{lemxi-eta}
If $f:\Gamma(q,\xi)\to\Gamma(q,\eta)$ is a graph isomorphism satisfying $\pi_\eta f=\pi_\xi$, where $\pi_\xi(H_\xi x)=Sx$ and $\pi_\eta(H_\eta x)=Sx$,
then $\xi=\eta$.
\end{lemma}
\begin{proof}
Let $R_\xi(h)$ and $R_\eta(h)$ denote right multiplication by $h\in G$ on the vertices of respective covers $\Gamma(q,\xi)$ and $\Gamma(q,\eta)$. By the assumption, the automorphisms $fR_\xi(h)f^{-1}$ and $R_\eta(h)$ induce the same permutation on $\Delta_n$. Then $R_\eta(h)^{-1}fR_\xi(h)f^{-1}=\alpha(h)\in \langle \tau_\eta\rangle$. Note that $\alpha: G\mapsto \langle\tau_{\eta}\rangle$ is a homomorphism (since $\tau_{\eta}$ is centralized by $G$). But $G$ is
simple, so $\alpha$ must be trivial.

As   $f(H_\xi)=H_\eta s$ for some $s\in S$,  equality of the stabilizers of these vertices for the two right $G$-actions gives $$ H_\xi=s^{-1}H_\eta s=H_\eta,$$ because $S$ is abelian. Therefore the two nonzero binary functionals $\ell_\xi,\ell_\eta$ coincide, so $\xi=\eta$.
 \end{proof}

\noindent\textit{Proof of Theorems~\ref{thm_full_aut} and~\ref{thm_iso}.}
By Lemma~\ref{lem:deck-intrinsic}, each automorphism of $\Gamma(\xi,q)$ preserves the  partition into $\tau_\xi$-opposite pairs and hence induces an automorphism of $\Delta_n$. Clearly, the kernel of the induced action  is exactly $\langle\tau_\xi\rangle$. By~\eqref{eq:aut-quotient},
every induced automorphism is a product of an element of $G$ and a power of a Frobenius field automorphism.

Note that $\pi_\xi R_\xi(g)=R_\Delta(g)\pi_\xi$ for all $g\in G$. Then, entrywise Frobenius  automorphism gives an isomorphism
$$ F_j:\Gamma(q,\xi)\longrightarrow\Gamma(q,\xi^{2^j}),  \qquad H_\xi x\longmapsto H_{\xi^{2^j}}x^{(2^j)},$$
such that $\pi_{\xi^{2^j}} F_j=\overline F_j \pi_\xi$, where $$\overline F_j:\Delta_n\longrightarrow\Delta_n, \qquad Sx\longmapsto Sx^{(2^j)}.$$
is the automorphism of $\Delta_n$ corresponding to $F_j$.
If that Frobenius power also lifts to $\Aut(\Gamma(q,\xi))$, composition with $F_j^{-1}$ gives an isomorphism  between $\Gamma(q,\xi^{2^j})$ and $\Gamma(q,\xi)$ that induces an identity automorphism of $\Delta_n$, so  Lemma~\ref{lemxi-eta} implies $\xi^{2^j}=\xi$.
Conversely, whenever $\xi^{2^j}=\xi$, the map $F_j$ is an automorphism of $\Gamma(q,\xi)$: it sends an admissible $c$ to
another admissible choice $c^{2^j}$ and the graph does not depend on that choice. Thus the image of $\Aut(\Gamma(q,\xi))$ in
$\Aut(\Delta_n)$ is exactly $G\rtimes\mathcal F_\xi$.

The lifts of $G$ and $\mathcal F_\xi$ form a subgroup which has trivial intersection with $\langle\tau_\xi\rangle$, and the involution $\tau_\xi$ commutes with both, since, for $\sigma\in\mathcal F_\xi$ one has $\ell_\xi(\sigma(c)+c)=0$, so $H_\xi u(\sigma(c))=H_\xi u(c)$.
Then, vertex and arc stabilizer orders follow by orbit-stabilizer theorem.
This proves Theorem~\ref{thm_full_aut}.

Let $f:\Gamma(q,\xi)\to\Gamma(q,\eta)$ be an isomorphism. Since the partitions into $\tau$-opposite pairs are  determined by the graph structure, $f$ maps the $\tau_\xi$-opposite pairs onto the $\tau_\eta$-opposite pairs. Hence $f$ induces an automorphism $\bar f\in\operatorname{Aut}(\Delta_n)$ such that
$$\pi_\eta f=\bar f\pi_\xi.$$

By the established description of $\operatorname{Aut}(\Delta_n)$, we may write
$$\bar f=R_\Delta(g)\varphi_j$$
for some $g\in G$ and $0\le j<n$, where
$$\varphi_j(Sx)=Sx^{(2^j)}.$$

Put
$$f_0=R_\eta(g)^{-1}f.$$

Since
$$\pi_\eta R_\eta(g)=R_\Delta(g)\pi_\eta,$$
we obtain
$$\pi_\eta f_0=R_\Delta(g)^{-1}\pi_\eta f=R_\Delta(g)^{-1}\bar f\pi_\xi=\varphi_j\pi_\xi.$$

On the other hand, the Frobenius isomorphism
$$F_j:\Gamma(q,\xi)\longrightarrow\Gamma(q,\xi^{2^j})$$
satisfies
$$\pi_{\xi^{2^j}}F_j=\varphi_j\pi_\xi.$$

Consequently, the isomorphism
$$h=f_0F_j^{-1}:\Gamma(q,\xi^{2^j})\longrightarrow\Gamma(q,\eta)$$
satisfies
$$\pi_\eta h=\pi_\eta f_0F_j^{-1}=\varphi_j\pi_\xi F_j^{-1}=\varphi_j\varphi_j^{-1}\pi_{\xi^{2^j}}=\pi_{\xi^{2^j}}.$$

Thus $h$ is an isomorphism of covers over the identity on $\Delta_n$. By Lemma~\ref{lemxi-eta},
$$\eta=\xi^{2^j}.$$

Conversely, if $\eta=\xi^{2^j}$, then $F_j$ is an isomorphism from $\Gamma(q,\xi)$ onto $\Gamma(q,\eta)$. This completes the proof of Theorem~\ref{thm_iso}.

\begin{corollary}
\label{cor:vertex-stabilizer}
For each vertex $v$ one has $$ \Aut(\Gamma(q,\xi))_v\simeq (Z_2)^n\rtimes Z_{n/d_\xi},$$
where the normal elementary abelian subgroup  acts regularly on the $q$ neighbours of $v$. In particular, none of these graphs is $2$-arc-transitive.
\end{corollary}
\begin{proof}
Put $T_\xi=\bigl\{\tau_\xi^{\ell_\xi (t)}R(u(t)):t\in\mathbb F_q\bigr\}$, where $R(u(t))\colon H_\xi x\longmapsto H_\xi x\,u(t)$. Note that, for all $t\in\F_q$, the permutations $\tau_\xi^{\ell_\xi(t)}R(u(t))$  fix the vertex $H_\xi$, form an elementary abelian group of order $q$, and act by translations on the neighbours, indexed by their lower rows $(1,z)$.
The Frobenius stabilizer $\mathcal F_\xi$ acts on this subgroup by $t\mapsto t^{2^j}$, and the orders in Theorem~\ref{thm_full_aut} show that these elements exhaust the vertex stabilizer. If the full group is $2$-arc-transitive, then it induces a $2$-transitive group on the $1$-neighbourhood of order at least $q(q-1)$. But $q(n/d_\xi)<q(q-1)$ for all $n\ge2$.
 \end{proof}

\begin{corollary}
\label{cor:orbit-counts}
For $0\le j<n$ let $e_j=\gcd(n,j)$ (so that $e_0=n$).
The numbers of abstract isomorphism classes in the trace-one
and nonzero trace-zero cases are respectively
\begin{equation}\label{eq:trace-one-iso-count}
 b_1(n)=\frac1n\sum_{\substack{0\le j<n\\n/e_j\ \mathrm{odd}}} 2^{e_j-1}
\end{equation}
and
\begin{equation}\label{eq:trace-zero-iso-count}
 b_0(n)=\frac1n\sum_{j=0}^{n-1}
 \begin{cases}
 2^{e_j-1}-1,&n/e_j\ \mathrm{odd},\\
 2^{e_j}-1,&n/e_j\ \mathrm{even}.
 \end{cases}
\end{equation}
\end{corollary}
\begin{proof}
By Theorem~\ref{thm_iso}, the isomorphism classes correspond to Frobenius orbits of $\xi$. For $\sigma_j:t\mapsto t^{2^j}$,
the fixed field is $\F_{2^{e_j}}$. The absolute trace of its members is $(n/e_j)$ times the trace from $\F_{2^{e_j}}$ to
$\F_2$. When $n/e_j$ is odd, exactly $2^{e_j-1}$ fixed coefficients have trace one, and $2^{e_j-1}-1$ nonzero fixed
coefficients have trace zero. When $n/e_j$ is even, none have trace one, and all $2^{e_j}-1$ nonzero fixed coefficients have
trace zero. Burnside's lemma now yields the required formulas.
 \end{proof}

 \section{Structural properties:  amply regularity, girth and spectrum}\label{sec:deza-girth}
 The results of this section hold for every $n\ge2$ and every nonzero
 trace coefficient $\xi\in\F_{2^n}^*$. We keep the notation
 $S,M,K,y,g$ fixed above and put $\delta_\xi=\Tr(\xi)$.

 \begin{lemma}\label{lem:cover-common-neighbours}
 Let $\pi:\widetilde\Lambda\to\Lambda$ be a graph covering, and suppose
 that any two distinct vertices of $\Lambda$ have at most one common
 neighbour. Then the same property holds in $\widetilde\Lambda$.
 Distinct vertices in a common fibre have no common neighbour. In
 particular, $\widetilde\Lambda$ has no $4$-cycles.
 \end{lemma}
 \begin{proof} Obvious.
 \end{proof}

 \begin{theorem}
 \label{thm:all-trace-deza}
 For every $q=2^n\ge4$ and $\xi\ne0$, the graph $\Gamma(q,\xi)$ is a
 connected, arc-transitive amply regular graph with parameters
 \begin{equation}\label{eq:all-trace-deza-parameters}
  (v,k,\lambda,\mu)=\bigl(2(q^2-1),\ q,\ \Tr(\xi),\ 1\bigr).
 \end{equation}

 \end{theorem}
 \begin{proof}
 The first two parameters and arc-transitivity follow from
 Propositions~\ref{prop:basic} and~\ref{prop:arc-transitive}; connectivity
 was established in Section~\ref{sec:parity}, with $q=4$ treated in
 Proposition~\ref{prop:q4-exception}. Since $\Ga$ covers $\Delta_n$,
combination of  Proposition~\ref{thm:involution} with  Lemma~\ref{lem:cover-common-neighbours}
  proves  $\mu=1$  in the covering graph.

 For edge regularity, take the base vertex $H_\xi$ and its neighbour
 with lower row $(1,D)$. Its binary invariant is $1+\ell_\xi(D)$.
 By Lemma~\ref{lem:two-step-param}, exactly one two-step word has that
 lower row and its invariant is $\ell_\xi(D)+\delta_\xi$. These invariants
 coincide precisely when $\delta_\xi=1$. Consequently the edge has one
 common neighbour in that case and none when $\delta_\xi=0$.
 Arc-transitivity gives the result for every edge.

\end{proof}

 \begin{theorem}\label{thm:all-trace-girth}
 For every $q=2^n\ge4$ and every $\xi\ne0$,
 \begin{equation}\label{eq:all-trace-girth}
 \operatorname{girth}\Gamma(q,\xi)=
 \begin{cases}
 3,&\Tr(\xi)=1,\\
 5,&\Tr(\xi)=0.
 \end{cases}
 \end{equation}
 \end{theorem}
 \begin{proof}
Every edge of $\Delta_n$  lies in a unique triangle. If $\delta_\xi=1$,
 Theorem~\ref{thm:all-trace-deza} shows that every edge
 of the cover lies in a unique triangle as well, so its girth is $3$.

 Suppose $\delta_\xi=0$. The graph is triangle-free by
 \eqref{eq:all-trace-deza-parameters}, and has no $4$-cycles by
 Lemma~\ref{lem:cover-common-neighbours}. Thus its girth is at least $5$.
Let us prove existence of a $5$-cycle.

 First let $q\ge8$. Choose any $C\in\F_q-\{0,1\}$ and any
 $D\in\F_q$. By Lemma~\ref{lem:two-step-param} there is exactly one
 lift $z$ of the lower row $(C,D)$ reached in two steps from $x=H_\xi$.
 It is not adjacent to $x$, since the images of $x,z$ are not adjacent
 in $\Delta_n$. By Lemma~\ref{lem:three-step-kloosterman}, both lifts
 of this row, in particular this same $z$, are reached by a three-step
 walk from $x$. Since $x,z$ are nonadjacent and the graph has no loops,
 this three-step walk is a simple path.
 Write the two paths as
 $$ x,v,z\qquad\text{and}\qquad x,r,s,z. $$
 They are internally disjoint ($v=r$ would give the triangle $v,s,z$,
 and $v=s$ would give the triangle $x,r,v$). Contradiction.
 Their union is therefore the $5$-cycle $x,v,z,s,r,x$.

 For $q=4$, the only nonzero trace-zero coefficient is $\xi=1$.
 Write $\F_4=\F_2(\omega)$ with $\omega^2+\omega+1=0$ and choose
 $y=u(\omega)$, so $\ell_1(\omega)=1$. With the fixed involution
 $g=\begin{pmatrix}0&1\\1&0\end{pmatrix}$, put $a_\omega=yg$.
 The Cayley--Hamilton identity gives
 $$ a_\omega^2=\omega a_\omega+I,\qquad  a_\omega^4=a_\omega+\omega I,\qquad a_\omega^5=I. $$
 Thus $a_\omega$ has order $5$. Since $|H_1|=2$, the cosets
 $$ H_1,\ H_1a_\omega,\ H_1a_\omega^2,\ H_1a_\omega^3,\ H_1a_\omega^4$$
 are distinct and consecutive cosets, including the last and the first,
 are adjacent. They form a $5$-cycle. This completes the proof for
 all $n\ge2$.
\end{proof}

 \begin{corollary}
 \label{cor:canonical-deza-girth}
 For every $n\ge2$, the canonical graph $\Gamma_n=\Gamma(2^n,1)$ is a
 Deza graph with parameters
 $$ \bigl(2(2^{2n}-1),\ 2^n,\ 1,\ 0\bigr),$$
 and
 \begin{equation}\label{eq:canonical-girth}
 \operatorname{girth}(\Gamma_n)=
 \begin{cases}3,&n\text{ odd},\\5,&n\text{ even}.\end{cases}
 \end{equation}
 Its diameter is $5$ for odd $n\ge3$ and $3$ for even $n\ge2$.
 \end{corollary}
 \begin{proof}
 Use $\Tr(1)=n\pmod2$ and Theorems~\ref{thm:all-trace-deza},
 \ref{thm:all-trace-girth}, \ref{thm:even-diameter}, and
 \ref{thm:odd-diameter}, together with
 Proposition~\ref{prop:q4-exception}.
\end{proof}

\begin{proposition}\label{prop:equitable}
Let $A=A_{\Ga}$, and let $P$ be the characteristic matrix of the
partition of $V\Ga$ into the $N=q^2-1$  $\tau_{\xi}$-opposite pairs. Then this partition is equitable with quotient
$\Delta_n$, and
$$AP=PA_{\Delta_n},\qquad P^{\mathsf T}P=2I, \qquad A_{\Delta_n}=\frac12P^{\mathsf T}AP.$$

After choosing a labeling $+$ and $-$ in each $\tau_\xi$-orbit, there
is a signed adjacency matrix $C_\xi$ of $\Delta_n$ such that
$$A \sim_{\mathrm{orth}} A_{\Delta_n}\oplus C_\xi.$$
Consequently,
$$\chi_{\Ga}(X)=\chi_{\Delta_n}(X)\det(XI-C_\xi).$$
\end{proposition}
\begin{proof}
The orbits of an automorphism form an equitable partition. In the present
case the quotient by $\langle\tau_\xi\rangle$ is precisely $\Delta_n$,
since
$$G/H_\xi\longrightarrow G/S$$
is the natural two-fold covering projection.

Since this is a graph covering, between two adjacent fibres there is
a perfect matching, while between nonadjacent fibres there are no edges.
Hence the quotient matrix is exactly $A_{\Delta_n}$. Therefore
$$AP=PA_{\Delta_n}.$$
Since every fibre has size $2$,
$$P^{\mathsf T}P=2I,$$
and multiplying $AP=PA_{\Delta_n}$ on the left by
$\frac12P^{\mathsf T}$ gives
$$A_{\Delta_n}=\frac12P^{\mathsf T}AP.$$

Choose in each fibre an ordering $(+,-)$. Relative to this ordering,
the adjacency matrix has the block form
$$A=
\begin{pmatrix}
A_0&A_1\\
A_1&A_0
\end{pmatrix},
$$
where $$A_0+A_1=A_{\Delta_n}.$$
Indeed, for each edge of $\Delta_n$, the corresponding perfect matching
is either parallel or crossed with respect to the chosen labels.

Put
$$C_\xi=A_0-A_1.$$
Then $C_\xi$ is a signed adjacency matrix of $\Delta_n$: its
$(u,v)$-entry is $0$ when $u$ and $v$ are nonadjacent in $\Delta_n$,
and is $\pm1$ when they are adjacent.

Now conjugate $A$ by the orthogonal matrix
$$Q=\frac1{\sqrt2}
\begin{pmatrix}
I&I\\
I&-I
\end{pmatrix}.
$$
A direct calculation gives
$$Q^{\mathsf T}AQ=
\begin{pmatrix}
A_0+A_1&0\\
0&A_0-A_1
\end{pmatrix}
=
A_{\Delta_n}\oplus C_\xi.
$$
Hence
$$A\sim_{\mathrm{orth}}A_{\Delta_n}\oplus C_\xi,$$
and therefore
$$\chi_{\Ga}(X)=\chi_{\Delta_n}(X)\det(XI-C_\xi).$$

Changing the labels $+$ and $-$ in an individual fibre conjugates
$C_\xi$ by a diagonal matrix with entries $\pm1$. Thus different
choices of fibre labels produce switching-equivalent signed adjacency
matrices and, in particular, do not change the spectrum of $C_\xi$.
\end{proof}

 \begin{proposition}\label{prop:quotient-spectrum}
 The spectrum of $\Delta_n$ is
 $$ q^1,\quad(-1)^q,\quad(\sqrt q)^{m_q},\quad(-\sqrt q)^{m_q},\qquad m_q=\frac{(q+1)(q-2)}2. $$
 All non-trivial eigenvalues of the quotient satisfy the Ramanujan bound.
 Thus, for every $\xi\ne0$,
 \begin{equation}\label{eq:signed-Ramanujan}
 \Ga\text{ is Ramanujan} \quad\Longleftrightarrow\quad \|C_\xi\|=\max_{\lambda\in\operatorname{Spec}(C_\xi)} |\lambda|\le2\sqrt{q-1}.
 \end{equation}
 \end{proposition}
 \begin{proof}
 The intersection matrix from Proposition~\ref{thm:involution} gives the
 four eigenvalues and their multiplicities. For $q\ge4$, we have
 $\sqrt q\le2\sqrt{q-1}$. Note that the graphs are not bipartite.
 Indeed, trace-one members contain triangles, and
 trace-zero members  have  even- and odd-distance
 part sizes $q^2-q+1$ and $q^2+q-3$, which differ for $q\ge4$ and hence cannot form bipartite parts.
 Thus the new part contains no exceptional trivial value $-q$.
\end{proof}


 It is worth mentioning that besides the trace-defined covers, $\Delta_n$ has its standard
 \emph{canonical bipartite double cover} $\widehat\Delta_n$ (of diameter $4$), also called its Kronecker cover;
 see \cite[Section~2]{ImrichPisanski2008}.

 \begin{remark}
 \label{rem:known-bipartite-spectrum}
 Put $A_\Delta=A_{\Delta_n}$. Then
 $$ A_{\widehat\Delta_n}=
 \begin{pmatrix}0&A_\Delta\\ A_\Delta&0\end{pmatrix}
 \simeq A_\Delta\oplus(-A_\Delta).
 $$
 Thus Proposition~\ref{prop:quotient-spectrum} gives
 \begin{equation}\label{eq:known-bipartite-spectrum}
 \Spec(\widehat\Delta_n)=\bigl\{q^1,(-q)^1,1^q,(-1)^q,(\sqrt q)^{2m_q},(-\sqrt q)^{2m_q}\bigr\},\qquad m_q=\frac{(q+1)(q-2)}2.
 \end{equation}
 After excluding the two trivial bipartite eigenvalues $\pm q$, the maximum
 absolute eigenvalue is $\sqrt q\le2\sqrt{q-1}$. The standard doubles are
 therefore Ramanujan by an immediate known-spectrum argument.
 \end{remark}
 \section*{Acknowledgements}
The research was conducted under the state assignment of the Institute of Mathematics and Mechanics of the Ural Branch of the Russian Academy of Sciences (project No. FUMF-2022-0003).

\end{document}